\documentclass[conference]{IEEEtran} 

\IEEEoverridecommandlockouts                              

\usepackage[margin=0.30cm]{caption}
\usepackage{epsf}
\usepackage{epsfig}
\usepackage{times}
\usepackage{amssymb}
\usepackage{amsmath}
\usepackage{amsfonts}
\usepackage{latexsym}
\usepackage{url}
\usepackage{mathrsfs}
\usepackage{algorithm}
\usepackage{caption}
\usepackage{datetime}
\usepackage{bm}
\usepackage{algorithm,algorithmic}
\usepackage{subfigure}
\usepackage{framed} 
\usepackage[framed]{ntheorem}
\usepackage{cite}
\usepackage{lipsum}
\usepackage{amstext}  
\usepackage{array}      
\usepackage{mathabx}
\usepackage{sansmath}
\usepackage{mathtools}
\usepackage{color}
\usepackage[utf8]{inputenc}
\usepackage{nomencl}
\usepackage{multirow}
\usepackage{tabu}
\usepackage{stfloats}

\makenomenclature

\newframedtheorem{frm-definition}{Definition}

\newtheorem{definition}{Definition}

\newtheorem{assumption}{Assumption}

\newcommand{\cL}{{\cal{L}}}
\newcommand{\cN}{{\cal N}}

\newcommand{\cE}{{\cal E}}
\newcommand{\cK}{{\cal K}}

\newcommand{\cY}{{\cal Y}}
\newcommand{\cO}{{\cal O}}

\allowdisplaybreaks[4]

\begin{document}
\title{Multi-Level Optimal Power Flow Solver in Large Distribution Networks}

\author{Xinyang Zhou$^{*}$, Yue Chen$^{*}$, Zhiyuan Liu$^{\dagger}$, Changhong Zhao$^{\ddagger}$, and Lijun Chen$^{\dagger}$\\
$^{*}$Power System Engineering Center, National Renewable Energy Laboratory, Golden, USA\\
$^{\dagger}$College of Engineering and Applied Science, University of Colorado, Boulder, USA\\
$^{\ddagger}$Department of Information Engineering, the Chinese University of Hong Kong, HK\\
Emails: \emph{\{xinyang.zhou, yue.chen\}@nrel.gov}, \emph{\{zhiyuan.liu, lijun.chen\}@colorado.edu}, \emph{chzhao@ie.cuhk.edu.hk}
\thanks{This work serves as an extension of our previous work \cite{zhou2019accelerated} published in \emph{IEEE Transactions on Power Systems}.}
\thanks{This work was authored in part by the National Renewable Energy Laboratory, operated by Alliance for Sustainable Energy, LLC, for the U.S. Department of Energy (DOE) under Contract No. DE-EE-0007998. Funding provided by U.S. Department of Energy Office of Energy Efficiency and Renewable Energy Solar Energy Technologies Office. The views expressed in the article do not necessarily represent the views of the DOE or the U.S. Government. The U.S. Government retains and the publisher, by accepting the article for publication, acknowledges that the U.S. Government retains a nonexclusive, paid-up, irrevocable, worldwide license to publish or reproduce the published form of this work, or allow others to do so, for U.S. Government purposes.}
}
\maketitle

\begin{abstract}
Solving optimal power flow (OPF) problems for large distribution networks incurs high computational complexity. We consider a large multi-phase distribution network of tree topology with a deep penetration of active devices. We divide the network into collaborating areas featuring subtree topology and subareas featuring subsubtree topology. We design a multi-level implementation of the primal-dual gradient algorithm to solve the voltage regulation OPF problems while preserving nodal voltage information and topological information within areas and subareas. Numerical results on a 4,521-node system verify that the proposed algorithm can significantly improve the computational speed without compromising any optimality. 
\end{abstract}


\section{Introduction}
As the renewable energy generation and smart electricity devices deepen their penetration in power systems, the modem distribution networks are becoming increasingly active. This trend not only suggests that we could exploit more flexibility in dispatchable devices to improve the overall system operation, but also brings challenges to architecting faster and more efficient distribution network operation paradigms that can take advantage of these active devices. Optimal power flow (OPF) is a powerful tool that can simultaneously determine the best operating points for all the dispatchable devices in the system under the constraints of operational limits by solving an optimization problem. 


The computational complexity of solving OPF, however, grows considerably when the size of the system increases. For example, the gradient algorithm for solving a voltage regulation problem (e.g., Eqs.~\eqref{eq:opt_phi}) requires computational complexity that is proportional to $N^2$, $N$ being the node number of the network. This means that, solving OPF for a 4,000-node system may take approximately 1,600 times longer than solving the same problem for a 100-node system. Other algorithms could experience similar scalability issues when used for large networks. This could lead to a slow response to fast-changing future distribution networks, especially in large systems.

To address this issue, various distributed algorithms have been developed to solve OPF in a more computationally efficient manner. The authors of \cite{dall2018optimala, tang2017real, hauswirth2017online, zhou2017discrete} use a central controller that is in charge of gathering, processing, and sending network-related information to coordinate distributed devices. The work in \cite{peng2016distributed, vsulc2014optimal,magnusson2019voltage, wu2018smart, kraning2014dynamic} does not need a central controller but require communication and computation capability among all neighbors. Also, a recent comprehensive survey on distributed optimization and control algorithms for electric power systems \cite{molzahn2017survey} contains more references.


Another promising line of work to enhance OPF solving efficiency proposes to divide a big network into smaller ones and thus solve smaller parts of the OPF problems in collaboration instead of the original large ones \cite{kroposki2018autonomous, conejo1998multi,nogales2003decomposition,lai2014decentralized}. The application of such a multi-area solution method in distribution networks, however, has been limited. Recent related works include \cite{zhang2018dynamic, chang2019saddle}. This work will continue to explore this method with novel algorithm that improves the computational efficiency without compromising optimality, which will be validated numerically on large systems.

In this work, we focus on a distribution network featuring tree topology, which is known to have fractal properties: any node within a tree together with all its children nodes makes a subtree that also features tree topology and thus inherits the same properties from the original tree; such properties are then passed on to subsubtrees within subtrees, and so on. This observation motivates us to develop a framework to exploit such a fractal pattern to improve the efficiency of solving OPF in large distribution networks. 

Following our previous work \cite{zhou2019accelerated}, this work continues to explore the fractal properties of distribution networks featuring tree topology, reveals the corresponding patterns in the sensitivity matrices between voltage magnitudes and nodal power injections in a multi-phases system, and proposes a multi-level solution method to solve optimal voltage regulation problems.  By design, the proposed algorithm can significantly improve the efficiency of solving large, convex OPF problems without losing optimality. Such highly efficient design is crucial for tracking the optimal operating setpoints in real time and for fast recovery from a blackout for large distribution systems with a large number of control nodes. We illustrate the numerical performance of the proposed multi-level algorithms on the primary side of a large distribution system based on the IEEE 8,500-node test feeder and the EPRI CKT7 test feeder. The results show identical convergence dynamics to the same optimal point from centralized algorithm, bi-level implementation, and tri-level implementation, but notably computational speed improvement by using multi-level algorithms.

The rest of this paper is organized as follows. Section~\ref{sec:model} models the distribution system, formulates an optimal voltage regulation problem, and introduces the primal-dual gradient algorithm for solving this OPF problem. Section~\ref{sec:multi} divides the large distribution system into subtrees and subsubtrees and proposes an equivalent multi-level implementation of the gradient algorithm. Section~\ref{sec:numerical} illustrates the numerical performance of the proposed algorithm, and Section~\ref{sec:conclusion} concludes this paper.

%

\section{System Model and OPF Solving}\label{sec:model}
\subsection{Network Modeling}
We consider a radial multi-phase distribution system denoted by $\mathcal{G}=\{\cN\cup\{0\}, \cE\}$ with the set $\cN$ collecting all $n$ buses excluding the substation bus 0, which is modeled as a slack bus, and the set $\cE$ collecting their connecting lines. We denote by $\cE_i$ the set of lines describing the unique path from the substation to bus $i$. Denote by $\mathfrak{i}:=\sqrt{-1}$ the unit imaginary number. Let $a,b,c$ denote the three phases, and 
$\Phi_i$ denote the set of phase(s) of bus $i\in\cN$. Obviously, $\Phi_i\subseteq \Phi_0=\{a,b,c\}$ of the substation bus 0 for all $i\in\cN$.
Define a subset $\cN^{\phi}\subseteq\cN$ collecting buses that have phase $\phi$. Denote by $p_i^{\phi}$, $q_i^{\phi}$, $V_i^{\phi}$ and $v_i^{\phi}$ the real power injection, the reactive power injection, the complex voltage phasor, and the squared voltage magnitude, respectively, of phase $\phi\in\Phi_i$ at bus $i\in\cN$. 
We make the following assumptions for a linearized power flow model.

\begin{assumption}\label{ass:3phase}
	The line losses are small and ignored, the magnitudes of the three-phase voltages are approximately equal, and the phase differences among the three-phase voltages are close to $2\pi/3$, i.e., $\frac{V_i^a}{V_i^b}\approx\frac{V_i^b}{V_i^c}\approx\frac{V_i^c}{V_i^a}\approx e^{\mathfrak{i} 2\pi/3}$.
\end{assumption}

\begin{figure}
	\centering
	\includegraphics[scale=0.38]{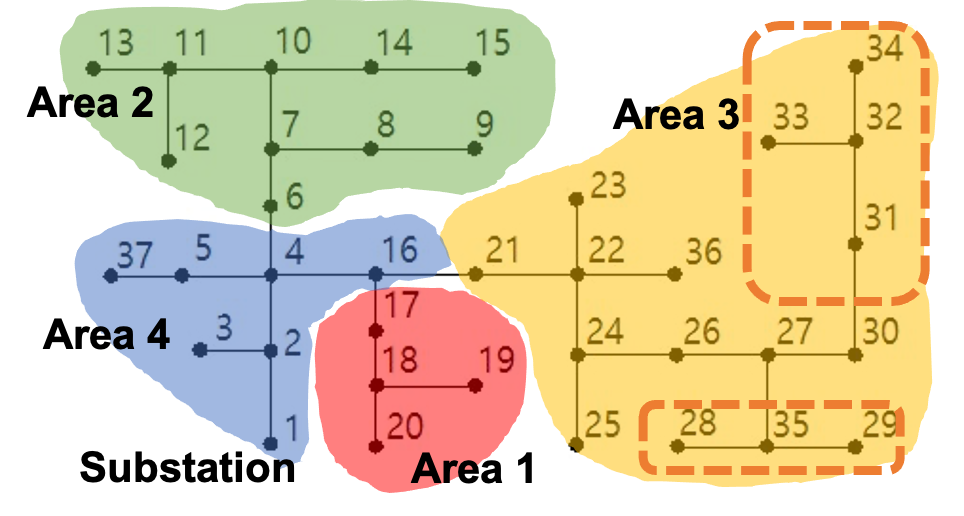}
	\caption{ Illustration of the subtree- and subsubtree-featuring structure based on IEEE 37-bus test feeder. Areas 1–3 ($\cN_1$--$\cN_3$) have subtree topology with their respective root buses $n_1^0=17$, $n_2^0=6$, and $n_3^0=21$. Area 4 is the remaining area with unclustered nodes ($\cN_0$). The circled subareas in Areas~3 are two subsubtrees ($\cN_{3_1}$ and $\cN_{3_2}$).}\label{fig:IEEE37}
\end{figure}

Denote by $z_{\zeta\xi}$ the impedance of line $(\zeta,\xi)\in\cE$. $z_{\zeta\xi}$ is a complex number for a single-phase line, a $2\times 2$ complex matrix for a two-phase line, or a $3\times 3$ complex matrix for a three-phase line. 
Denote by:
\begin{eqnarray}
 Z^{\varphi\phi}_{ij}=\sum_{(\zeta,\xi)\in\cE_i \cap \cE_j}z^{\varphi\phi}_{\zeta\xi}\in\mathbb{C}\nonumber
\end{eqnarray}
the summarized impedance (if $\varphi=\phi$) or the summarized mutual impedance (if $\varphi\neq\phi$) of the unique common path of buses~$i$ and bus~$j$ leading back to bus 0, and $\overline{Z}^{\varphi\phi}_{ij}$ its conjugate. For example, in Fig.~\ref{fig:IEEE37}, the common path of bus~10 and bus~27 leading back to bus~0 is line~(1,2) and line~(2,4).
Define:
\begin{eqnarray}
\bm{v}&:=&[[v_1^{\phi}]^{\top}_{\phi\in\Phi_1}, \ldots, [v_n^{\phi}]^{\top}_{\phi\in\Phi_N}]^{\top}\in\mathbb{R}^N\nonumber\\
\bm{p}&:=&[[p_1^{\phi}]_{\phi\in\Phi_1}^{\top}, \ldots, [p_n^{\phi}]^{\top}_{\phi\in\Phi_N}]^{\top}\in\mathbb{R}^N\nonumber\\
\bm{q}&:=&[[q_1^{\phi}]^{\top}_{\phi\in\Phi_1}, \ldots, [q_n^{\phi}]^{\top}_{\phi\in\Phi_N}]^{\top}\in\mathbb{R}^N\nonumber
\end{eqnarray}
the vectors of the squared voltage magnitudes, the real power injection, and the reactive power injection, with $N =\sum_{i\in\cN}|\Phi_i|$ adding up all phases at all buses. Here, $|\Phi_i|$ is the cardinality of set $\Phi_i$. 
To construct a convex and computationally trackable OPF problem, we use the following \emph{linearized} power flow equation \cite{gan2014convex, gan2016online} based on Assumption~\ref{ass:3phase}:
\begin{eqnarray}
\bm{v}=\mathbf{R}\bm{p}+\mathbf{X}\bm{q}+\tilde{\bm{v}},\label{eq:lindistflow_phi}
\end{eqnarray}
with a constant vector $\tilde{\bm{v}}\in\mathbb{R}^N$ and sensitivity matrices $\mathbf{R}, \mathbf{X}\in\mathbb{R}^{N\times N}$ respectively comprising the elements of:
\begin{eqnarray}
&&\partial_{p_j^{\phi}}v_i^{\varphi} =    2\mathfrak{Re} \big\{ \overline{Z}^{\varphi\phi}_{ij} \cdot \omega^{\varphi-\phi}\big\}\text{ and}\nonumber\\
&&\partial_{q_j^{\phi}}v_i^{\varphi} = -2 \mathfrak{Im} \big\{ \overline{Z}^{\varphi\phi}_{ij} \cdot \omega^{\varphi-\phi}\big\}.\nonumber
\end{eqnarray}
Here, $\omega=e^{-\mathfrak{i} 2\pi/3}$, $a,b,c=0,1,2$ when calculating $\varphi-\phi$, and
$\mathfrak{Re}\{\cdot\}$ and $\mathfrak{Im}\{\cdot\}$ are the real and imaginary parts of a complex number, respectively. 


\subsection{Problem Formulation and Gradient Algorithm}
To determine the optimal power setpoints of dispatchable devices while strictly constraining voltages to within specified bounds, we formulate a voltage regulation OPF problem for the multi-phase distribution system as follows:
\begin{subequations}\label{eq:opt_phi}
 	\begin{eqnarray}
 	\hspace{-5mm}	&\underset{\bm{p}_,\bm{q}}{\min} & \sum_{i\in\cN}\sum_{\phi\in\Phi_i}C_i^{\phi}(p_i^{\phi},q_i^{\phi}),\\   
 	\hspace{-5mm}	& \text{s.t.}&
 		\underline{\bm{v}} \leq \bm{v}(\bm{p},\bm{q}) \leq \overline{\bm{v}},\label{eq:voltreg_phi}\\
 	\hspace{-5mm}	&     & (p_i^{\phi},q_i^{\phi})\in\cY_i^{\phi},\phi\in\Phi_i, \forall i\in\cN, \label{eq:X_phi}
 	\end{eqnarray}
 \end{subequations}
where $C_i^{\phi}$ is a jointly strongly convex function in both $p_i^{\phi}$ and $q_i^{\phi}$ for phase $\phi\in\Phi$ of bus $i$, $\cY_i^{\phi}$ is its convex and compact feasible set, 
$\bm{v}(\bm{p},\bm{q})$ represents Eq.~\eqref{eq:lindistflow_phi}, and $\underline{\bm{v}}=\{\underline{v}_i^{\phi}\}_{i\in\cN}^{\phi\in\Phi_i}$ and  $\overline{\bm{v}}=\{\overline{v}_i^{\phi}\}_{i\in\cN}^{\phi\in\Phi_i}$ are the voltage bounds vectors. Associate dual variables $\underline{\bm{\mu}}=\{\underline{\mu}_i^{\phi}\}_{i\in\cN}^{\phi\in\Phi_i}$ and $\overline{\bm{\mu}}=\{\overline{\mu}_i^{\phi}\}_{i\in\cN}^{\phi\in\Phi_i}$ with the left and right sides of the constraints \eqref{eq:voltreg_phi}, respectively, and we write the regularized Lagrangian of \eqref{eq:opt_phi} as:
\begin{align}
\cL(\bm{p},\bm{q};\overline{\bm{\mu}},\underline{\bm{\mu}})=&\sum_{i\in\cN}\sum_{\phi\in\Phi_i}C_i^{\phi}(p_i^{\phi},q_i^{\phi})+\underline{\bm{\mu}}^{\top}(\underline{\bm{v}}-\bm{v}(\bm{p},\bm{q})) \nonumber\\
&\hspace{4mm}+\overline{\bm{\mu}}^{\top}(\bm{v}(\bm{p},\bm{q})-\overline{\bm{v}})-\frac{\eta}{2}\|\bm{\mu}\|^2_2,
\label{eq:langr2}
\end{align}
where a regularization term $-\frac{\eta}{2}\|\bm{\mu}\|^2_2$ with a small $\eta>0$ is added to improve convergence properties. A minor bounded discrepancy between the regularized Lagrangian and the original Lagrangian is introduced because of the regularization term \cite{koshal2011multiuser}. 

We implement the primal-dual gradient algorithm for solving the unique saddle point of \eqref{eq:langr2} with a stepsize $\epsilon$ as:
\begin{subequations}\label{eq:primaldual3}
\begin{eqnarray}
\hspace{-6mm}p_i^{\phi}(t+1)\hspace{-2mm}&=&\hspace{-2mm}\Big[p^{\phi}_i(t)-\epsilon\Big(\partial_{p^{\phi}_i} C^{\phi}_i(p^{\phi}_i(t),q^{\phi}_i(t))+\nonumber\\
\hspace{-6mm}&&\hspace{2mm}\sum_{j\in\cN}\sum_{\varphi\in\Phi_j}	\partial_{p_i^{\phi}}v_j^{\varphi}\big(\overline{\mu}^{\varphi}_j(t)-\underline{\mu}^{\varphi}_j(t)\big)\Big) \Big]_{\cY^{\phi}_i},\label{eq:mpPi}\\
\hspace{-6mm}q_i^{\phi}(t+1)\hspace{-2mm}&=&\hspace{-2mm}\Big[q^{\phi}_i(t)-\epsilon\Big(\partial_{q^{\phi}_i} C^{\phi}_i(p^{\phi}_i(t),q^{\phi}_i(t))+\nonumber\\
\hspace{-6mm}&&\hspace{2mm}\sum_{j\in\cN}\sum_{\varphi\in\Phi_j}	\partial_{q_i^{\phi}}v_j^{\varphi}\big(\overline{\mu}^{\varphi}_j(t)-\underline{\mu}^{\varphi}_j(t)\big)\Big) \Big]_{\cY^{\phi}_i},\label{eq:mpQi}\\
\hspace{-6mm}\underline{\mu}^{\phi}_i(t+1)\hspace{-2mm}&=&\hspace{-2mm}[\underline{\mu}^{\phi}_i(t)+\epsilon (\underline{v}^{\phi}_i-v^{\phi}_i(t)-\eta \underline{\mu}^{\phi}_i(t))]_+,\label{eq:mpmu1}\\
\hspace{-6mm}\overline{\mu}^{\phi}_i(t+1)\hspace{-2mm}&=&\hspace{-2mm}[\overline{\mu}^{\phi}_i(t)+\epsilon (v^{\phi}_i(t)-\overline{v}^{\phi}_i-\eta \overline{\mu}^{\phi}_i(t))]_+,\label{eq:mpmu2}\\
\hspace{-6mm}\bm{v}(t+1)\hspace{-2mm}&=&\hspace{-2mm} \mathbf{R}\bm{p}(t+1)+\mathbf{X}\bm{q}(t+1)+\tilde{\bm{v}},\label{eq:linearphi}
\end{eqnarray}
\end{subequations}
where \eqref{eq:mpPi}--\eqref{eq:mpmu2} are for all $\phi\in\Phi_i$ and all $i\in\cN$. Proof of \eqref{eq:primaldual3} converging asymptotically to its unique saddle point of Lagrangian~\eqref{eq:langr2} given a small enough stepsize $\epsilon$ can be found in numerous references (e.g., \cite{bertsekas1989parallel,zhou2019accelerated}) and is omitted here.

\subsection{Feedback from Nonlinear Power Flow}
What has been presented so far is based on the linearized power flow equation~\eqref{eq:lindistflow_phi} assuming a lossless and balanced three-phase system. As will be shown in Section~\ref{sec:multi}, such linearized model is crucial to obtain an important structure in the sensitivity matrices, allowing us to design a hierarchical algorithm.
However, unrealistic assumptions like this will inevitably cause modeling errors in practice. For example, when we obtain an optimal solution $(\bm{p}^*, \bm{q}^*)$ and the corresponding $\bm{v}^*$ through dynamics~\eqref{eq:primaldual3}, and implement the setpoints in real world featuring nonlinear power flow, the resultant voltage magnitudes may be different than $\bm{v}^*$, putting the actual system at risk of voltage violation.

To compensate for the linearization errors and to avoid voltage violation in practice, a model-based feedback approach can be used by replacing Eq.~\eqref{eq:linearphi} with its nonlinear counterpart so that the decision variables are updated based on voltages from nonlinear power flow model. We refer the detailed performance characterization to \cite{zhou2019accelerated}. This feedback approach is also applied in the numerical examples in Section~\ref{sec:numerical}, where nonlinear power flow is calculated each iteration.

\section{Multi-Level OPF Algorithm}\label{sec:multi}

\subsection{Motivation}
Implementing Eqs.~\eqref{eq:primaldual3} for multiple iterations until convergence calls for a great amount of computation, especially for large distribution networks. Particularly, the computational complexity for 
\begin{subequations}\label{eq:complex}
\begin{eqnarray}
&&\sum_{j\in\cN}\sum_{\varphi\in\Phi_j}	\partial_{p_i^{\phi}}v_j^{\varphi}\big(\overline{\mu}^{\varphi}_j(t)-\underline{\mu}^{\varphi}_j(t)\big)\Big) \Big]_{\cY^{\phi}_i}\\
&&\sum_{j\in\cN}\sum_{\varphi\in\Phi_j}	\partial_{q_i^{\phi}}v_j^{\varphi}\big(\overline{\mu}^{\varphi}_j(t)-\underline{\mu}^{\varphi}_j(t)\big)\Big) \Big]_{\cY^{\phi}_i}
\end{eqnarray}
\end{subequations}
in Eqs.~\eqref{eq:mpPi}--\eqref{eq:mpQi} for all $\phi\in\Phi_i$ and all $i\in\cN$ increases quadratically as the size of the network $N$ increases. Moreover, certain areas within a distribution network may want to preserve local nodal and topological information for security or privacy reasons.


For these reasons, we redesign the implementation of Eqs.~\eqref{eq:primaldual3} for radial distribution networks to significantly increase the computational efficiency and preserve privacy for areas with needs to protect nodal or topological information. Specifically, when we divide the large system geographically/administratively into areas featuring subtrees, we can take advantage of the network structure as well as the linearization model \eqref{eq:lindistflow_phi} to simplify the calculation of the computationally heavy coupling parts without compromising any performance based on the collaboration of divided areas. In this way, the computational efficiency is significantly improved while the nodal and topological information within areas are preserved from outside of the area. Next, we apply similar techniques to each area (subtree) to further boost the computational efficiency and preserve privacy within subareas (subsubtrees).


\subsection{Multi-Level Clustering}
Based on the pattern of the sensitivity matrix related to the tree/subtree structure in the distribution network, we group all nodes of distribution network $\cN$ into two categories: $K$ nonoverlapping subtrees indexed by $k\in\cK:=\{1,\ldots,K\}$ with $\cN_k\in\cN$ collecting all nodes within subtree $k$, and a set $\cN_0$ collecting all the other ``unclustered" nodes in $\cN$. Here, we have $\cup_{k\in\cK}\cN_k\cup\cN_0=\cN$ and $\cN_j \cap \cN_k=\emptyset,\forall j\neq k$. Denote the root node of subtree $\cN_k$ by $n_k^0$.

Notice that a subtree inherits all the properties from the original tree; therefore, we can divide one subtree into ``subsubtrees" for reasons that will become clear soon. We divide area $\cN_k$ into non-overlapping subsubtrees indexed by $k_m\in\cK_k := \{k_1,k_2,\ldots\}$ and their corresponding bus sets denoted by $\cN_{k_m}$ together with $\cN_{k_0}$ collecting the remaining unclustered buses.

\subsection{Design Intuition}
First, we formally define subtrees within a tree topology.
\begin{definition}
A subtree of a tree consists of a node in the tree, all this node's descendants, and the connecting lines among them.
\end{definition}
Because the sensitivity matrices $\mathbf{R}$ and $\mathbf{X}$ are built based on $Z_{ij}=\sum_{(\zeta,\xi)\in\cE_i \cap \cE_j}z^{\varphi\phi}_{\zeta\xi}$ that depends on the common path of bus $i$ and $j$ to the substation, i.e., $\cE_i \cap \cE_j$, and because two non-overlapping subtrees hold the same common path $\cE_i \cap \cE_j$ for any of their respective buses $i$ and bus $j$, we conclude that:
\begin{subequations}\label{eq:intuiiton}
	\begin{eqnarray}
	\hspace{-14mm}&&\partial_{p_j^{\phi}}v_i^{\varphi} = \partial_{p_{n_h^0}^{\phi}}v_{n_k^0}^{\varphi} =\hspace{2.5mm}2\mathfrak{Re} \big\{ \overline{Z}^{\varphi\phi}_{n_k^0 n_h^0} \cdot \omega^{\varphi-\phi}\big\}\label{eq:intuiiton_a}\\
	\hspace{-14mm}&&\partial_{q_j^{\phi}}v_i^{\varphi} = \partial_{q_{n_h^0}^{\phi}}v_{n_k^0}^{\varphi} =-2\mathfrak{Im} \big\{ \overline{Z}^{\varphi\phi}_{n_k^0 n_h^0} \cdot \omega^{\varphi-\phi}\big\}
	\end{eqnarray}
\end{subequations}
hold for any $i\in\cN_k$ and any $j\in\cN_h$, where $\cN_k$ and $\cN_h$ respectively denote the sets collecting buses of two non-overlapping subtrees indexed by $k, h\in\cK$ with their respective root buses indexed by $n_k^0$ and $n_h^0$, with set $\cK$ collecting all subtree indexes. See Fig.~\ref{fig:multilevel} for an illustration, where the white dashed line is the common path of any node in Area~1 and any node in Area 3---marked by blue triangles---to the substation. Because subtrees inherit all features from the original tree, we can apply similar results to them.

\subsection{Methods}

\subsubsection{Bi-Level Implementation}
Based on such intuition, we can equivalently rewrite the last term of \eqref{eq:mpPi} for clustered node $i\in\cN_k$ with phase $\phi\in\Phi_i$ as the following three parts:
\begin{subequations}\label{eq:bip}
	\begin{eqnarray}
	&&\hspace{-12mm}\sum_{j\in\cN}\sum_{\varphi\in\Phi_j}	\partial_{p_i^{\phi}}v_j^{\varphi}\big(\overline{\mu}^{\varphi}_j(t)-\underline{\mu}^{\varphi}_j(t)\big)\Big) \Big]_{\cY^{\phi}_i}=\label{eq:part0}\\
	&&\hspace{-12mm}2\mathfrak{Re} \Big\{ \sum_{\varphi\in\Phi_0}\omega^{\varphi-\phi}\hspace{-4mm}\sum_{j\in\cN^{\varphi}\cap\cN_k}\hspace{-4mm}  \overline{Z}^{\varphi\phi}_{ji} \big(\overline{\mu}^{\varphi}_j(t)-\underline{\mu}^{\varphi}_j(t)\big)\Big\}+\label{eq:parta}\\
	&&\hspace{-12mm}2\mathfrak{Re} \Big\{ \sum_{\varphi\in\Phi_0}\omega^{\varphi-\phi}\hspace{-4mm}\sum_{\substack{n_h^0 \in N^\varphi\\h\in\cK,h\neq k }}\hspace{-3mm}\overline{Z}^{\varphi\phi}_{n_h^0 n_k^0}\hspace{-4mm}\sum_{j\in\cN^{\varphi}\cap\cN_h}\hspace{-4mm}\big(\overline{\mu}^{\varphi}_j(t)-\underline{\mu}^{\varphi}_j(t)\big)\Big\}+\label{eq:partb}\\
	&&\hspace{-12mm}2\mathfrak{Re} \Big\{ \sum_{\varphi\in\Phi_0}\omega^{\varphi-\phi}\hspace{-4mm}\sum_{j\in\cN^{\varphi}\cap\cN_0}\hspace{-4mm}\overline{Z}^{\varphi\phi}_{j n_k^0}\big(\overline{\mu}^{\varphi}_j(t)-\underline{\mu}^{\varphi}_j(t)\big)\Big)\Big\},\label{eq:partc}
	\end{eqnarray}
\end{subequations}
where \eqref{eq:parta}, \eqref{eq:partb}, and \eqref{eq:partc}, are the components of \eqref{eq:part0} from nodes within subtree $\cN_k$, nodes within all other subtrees $\cN_h, h\neq k$, and nodes not in any subtree $\cN_0$, respectively.
Computation-wise, \eqref{eq:parta} and \eqref{eq:partc} are same as before; however, \eqref{eq:partb} significantly simplifies the original computation by reducing a large number of repetitive computations by applying the result in \eqref{eq:intuiiton_a}, and \eqref{eq:partc} captures the influence from buses outside any divided subtree area (e.g., Area 5 in Fig.~\ref{fig:multilevel}) whose buses are collected in set $\cN_0$ for completeness.

The last term in \eqref{eq:mpQi} can be processed similarly as:
\begin{subequations}\label{eq:biq}
	\begin{eqnarray}
	&&\hspace{-12mm}\sum_{j\in\cN}\sum_{\varphi\in\Phi_j}	\partial_{q_i^{\phi}}v_j^{\varphi}\big(\overline{\mu}^{\varphi}_j(t)-\underline{\mu}^{\varphi}_j(t)\big)\Big) \Big]_{\cY^{\phi}_i}=\label{eq:part0q}\\
	&&\hspace{-12mm}-2\mathfrak{Im} \Big\{ \sum_{\varphi\in\Phi_0}\omega^{\varphi-\phi}\hspace{-4mm}\sum_{j\in\cN^{\varphi}\cap\cN_k}\hspace{-4mm}  \overline{Z}^{\varphi\phi}_{ji} \big(\overline{\mu}^{\varphi}_j(t)-\underline{\mu}^{\varphi}_j(t)\big)\Big\}\label{eq:partaq}\\
	&&\hspace{-12mm}-2\mathfrak{Im} \Big\{ \sum_{\varphi\in\Phi_0}\omega^{\varphi-\phi}\hspace{-4mm}\sum_{\substack{n_h^0 \in N^\varphi\\h\in\cK,h\neq k }}\hspace{-3mm}\overline{Z}^{\varphi\phi}_{n_h^0 n_k^0}\hspace{-4mm}\sum_{j\in\cN^{\varphi}\cap\cN_h}\hspace{-4mm}\big(\overline{\mu}^{\varphi}_j(t)-\underline{\mu}^{\varphi}_j(t)\big)\Big\}\label{eq:partbq}\\
	&&\hspace{-12mm}-2\mathfrak{Im} \Big\{ \sum_{\varphi\in\Phi_0}\omega^{\varphi-\phi}\hspace{-4mm}\sum_{j\in\cN^{\varphi}\cap\cN_0}\hspace{-4mm}\overline{Z}^{\varphi\phi}_{j n_k^0}\big(\overline{\mu}^{\varphi}_j(t)-\underline{\mu}^{\varphi}_j(t)\big)\Big)\Big\}.\label{eq:partcq}
	\end{eqnarray}
\end{subequations}

\subsubsection{Tri-Level Implementation}
Consider \eqref{eq:parta} and notice that it has a form similar to the original coupling term. Hence, the computational complexity of $\cO(|\cN_k|^2)$ for \eqref{eq:parta} can be high for large $|\cN_k|$. Fortunately, we can apply similar tricks to divide area $\cN_k$ into subsubtrees indexed by $k_m\in\cK_k = \{k_1,k_2,\ldots\}$ and their corresponding buses sets denoted by $\cN_{k_m}$ together with $\cN_{k_0}$ collecting the remaining buses. We then equivalently rewrite \eqref{eq:parta} for $i\in\cN_{k_m}\subset\cN_k, m\in\cK_k, k\in\cK$ as:
\begin{subequations}\label{eq:trilevel}
	\begin{eqnarray}
	&&\hspace{-12mm} 2\mathfrak{Re} \Big\{ \sum_{\varphi\in\Phi_0}\omega^{\varphi-\phi}\hspace{-4mm}\sum_{j\in\cN^{\varphi}\cap\cN_k}\hspace{-4mm}  \overline{Z}^{\varphi\phi}_{ji} \big(\overline{\mu}^{\varphi}_j(t)-\underline{\mu}^{\varphi}_j(t)\big)\Big\}=\label{eq:parta0}\\
	&&\hspace{-12mm}2\mathfrak{Re} \Big\{ \sum_{\varphi\in\Phi_0}\omega^{\varphi-\phi}\hspace{-4mm}\sum_{j\in\cN^{\varphi}\cap\cN_{k_m}}\hspace{-4mm}  \overline{Z}^{\varphi\phi}_{ji} \big(\overline{\mu}^{\varphi}_j(t)-\underline{\mu}^{\varphi}_j(t)\big)\Big\}+\label{eq:parta1}\\
	&&\hspace{-12mm}2\mathfrak{Re} \Big\{ \sum_{\varphi\in\Phi_0}\omega^{\varphi-\phi}\hspace{-6mm}\sum_{\substack{n_{k_{m'}}^0\!\! \in \cN^\varphi\\k_m'\in\cK_k,m'\neq m }}\hspace{-6mm}\overline{Z}^{\varphi\phi}_{n_{k_{m'}}^0\!\! n_{k_m}^0}\hspace{-4mm}\sum_{j\in\cN^{\varphi}\cap\cN_{k_{m'}}}\hspace{-7mm}\big(\overline{\mu}^{\varphi}_j(t)-\underline{\mu}^{\varphi}_j(t)\big)\Big\}+\label{eq:parta2}\\
	&&\hspace{-12mm}2\mathfrak{Re} \Big\{ \sum_{\varphi\in\Phi_0}\omega^{\varphi-\phi}\hspace{-4mm}\sum_{j\in\cN^{\varphi}\cap\cN_{k_0}}\hspace{-4mm}\overline{Z}^{\varphi\phi}_{j n_{k_m}^0}\big(\overline{\mu}^{\varphi}_j(t)-\underline{\mu}^{\varphi}_j(t)\big)\Big)\Big\},\label{eq:parta3}
	\end{eqnarray}
\end{subequations}
where $n_{k_m}^0$ denotes the root bus of subsubtree $k_m$. Now \eqref{eq:parta2} further simplifies the computation for \eqref{eq:parta} similarly as \eqref{eq:partb} does for \eqref{eq:mpPi}. Note that the nodal information within subsubtrees is aggregated without being exposed.  

Then Eq.~\eqref{eq:parta} for $i\in\cN_{k_m}\subset\cN_k, m\in\cK_k, k\in\cK$ can be similarly calculated as:
 \begin{subequations}\label{eq:trilevelq}
	\begin{eqnarray}
	&&\hspace{-12mm} -2\mathfrak{Im} \Big\{ \sum_{\varphi\in\Phi_0}\omega^{\varphi-\phi}\hspace{-4mm}\sum_{j\in\cN^{\varphi}\cap\cN_k}\hspace{-4mm}  \overline{Z}^{\varphi\phi}_{ji} \big(\overline{\mu}^{\varphi}_j(t)-\underline{\mu}^{\varphi}_j(t)\big)\Big\}=\label{eq:parta0q}\\
	&&\hspace{-12mm}-2\mathfrak{Im} \Big\{ \sum_{\varphi\in\Phi_0}\omega^{\varphi-\phi}\hspace{-4mm}\sum_{j\in\cN^{\varphi}\cap\cN_{k_m}}\hspace{-4mm}  \overline{Z}^{\varphi\phi}_{ji} \big(\overline{\mu}^{\varphi}_j(t)-\underline{\mu}^{\varphi}_j(t)\big)\Big\}\label{eq:parta1q}\\
	&&\hspace{-12mm}-2\mathfrak{Im} \Big\{ \sum_{\varphi\in\Phi_0}\omega^{\varphi-\phi}\hspace{-6mm}\sum_{\substack{n_{k_{m'}}^0\!\! \in \cN^\varphi\\k_m'\in\cK_k,m'\neq m }}\hspace{-6mm}\overline{Z}^{\varphi\phi}_{n_{k_{m'}}^0\!\! n_{k_m}^0}\hspace{-4mm}\sum_{j\in\cN^{\varphi}\cap\cN_{k_{m'}}}\hspace{-7mm}\big(\overline{\mu}^{\varphi}_j(t)-\underline{\mu}^{\varphi}_j(t)\big)\Big\}\label{eq:parta2q}\\
	&&\hspace{-12mm}-2\mathfrak{Im} \Big\{ \sum_{\varphi\in\Phi_0}\omega^{\varphi-\phi}\hspace{-4mm}\sum_{j\in\cN^{\varphi}\cap\cN_{k_0}}\hspace{-4mm}\overline{Z}^{\varphi\phi}_{j n_{k_m}^0}\big(\overline{\mu}^{\varphi}_j(t)-\underline{\mu}^{\varphi}_j(t)\big)\Big)\Big\}.\label{eq:parta3q}
	\end{eqnarray}
\end{subequations}

Eqs.~\eqref{eq:bip}--\eqref{eq:trilevelq} provide a more efficient way to carry out the most computationally heavy terms in Eqs.~\eqref{eq:mpPi}--\eqref{eq:mpQi}. While the proposed method has reduced the computational complexity, its convergence dynamics and equilibrium point stay the same because the new method is mathematically equivalent to the original primal-dual gradient algorithm. Such design gives us ``free" acceleration without losing performance, even if it is not implemented in a parallel manner among areas and subareas, as will be shown in Section~\ref{sec:numerical}. Furthermore, considering the fractal properties of tree/subtree---i.e., subtrees have tree structures and thus can be further divided into subsubtrees where everything still holds---we can theoretically apply such layering structure as deep as we wish, if it is needed and if the network allows.

\subsection{Privacy Preservation}
In the original primal-dual gradient algorithm~\eqref{eq:primaldual3}, a central controller needs to know the dual variables for all buses and the topology of the entire network to compute the coupling parts~\eqref{eq:complex}; therefore, every bus must share its nodal voltage or local dual variables with the central controller every iteration, together with its topological information. Meanwhile, when computing \eqref{eq:complex} through Eqs.~\eqref{eq:bip}--\eqref{eq:biq}, all nodal and topological information can be preserved within each area, and only aggregated area information needs to be shared with the central controller. For example, for all buses within area $k$, they only need to share their aggregated dual variables $\sum_{j\in\cN^{\varphi}\cap\cN_k}\big(\overline{\mu}^{\varphi}_j(t)-\underline{\mu}^{\varphi}_j(t)\big)\Big\}$ for all phases $\varphi$ and the topological information related to its root buses $n_k^0$ to outside of area $k$. For buses within subareas, detailed information of subareas is similarly preserved within.

\begin{figure}[h]
	\centering
	\includegraphics[scale=0.2]{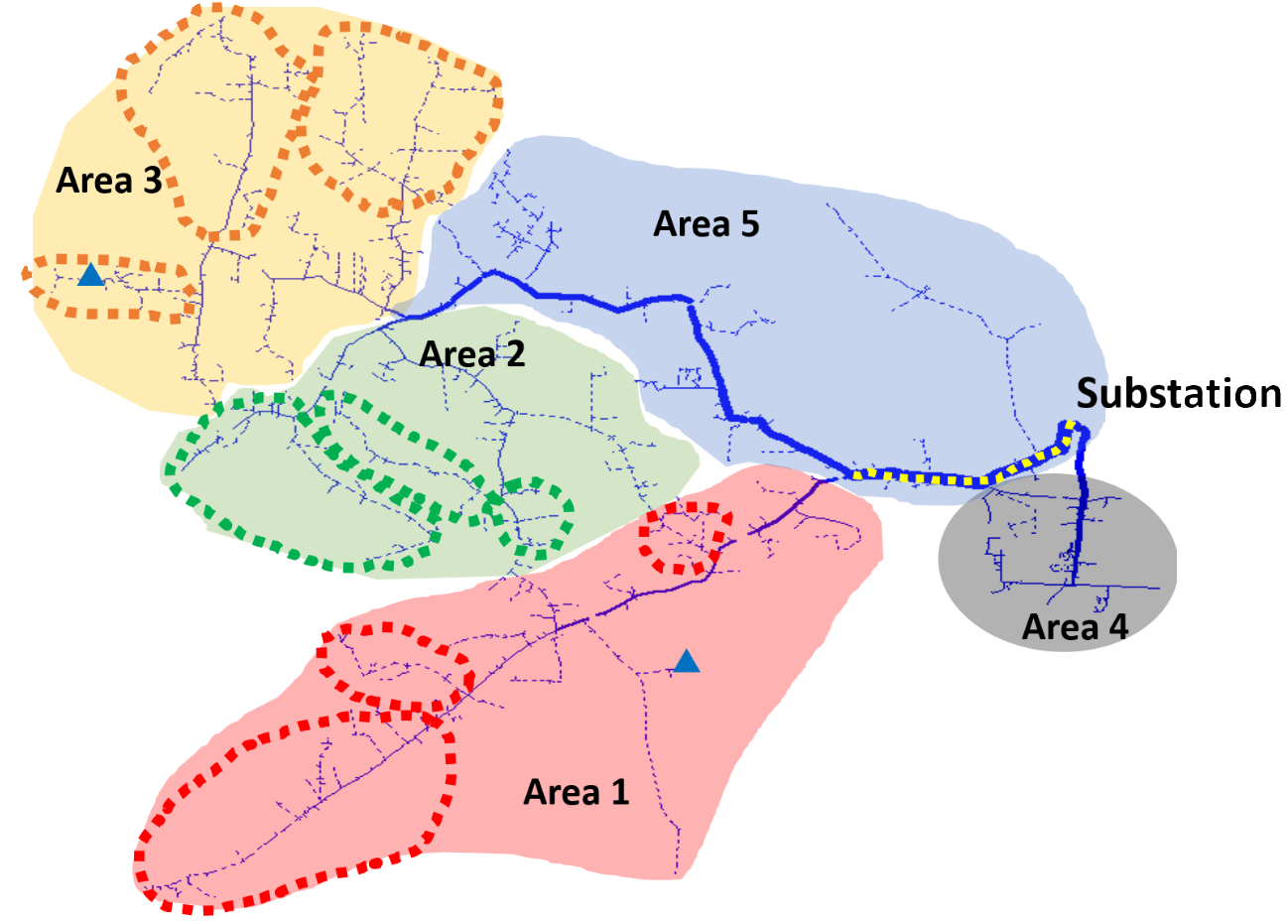}
	\caption{Areas and subareas (within dashed lines) of the 4,521-node network used in Section~\ref{sec:numerical}. Areas 1--3 have subtree topology, and all the subareas feature subsubtree topology.}\label{fig:multilevel}
\end{figure}

\begin{figure*}
	\centering
	\includegraphics[scale=0.375]{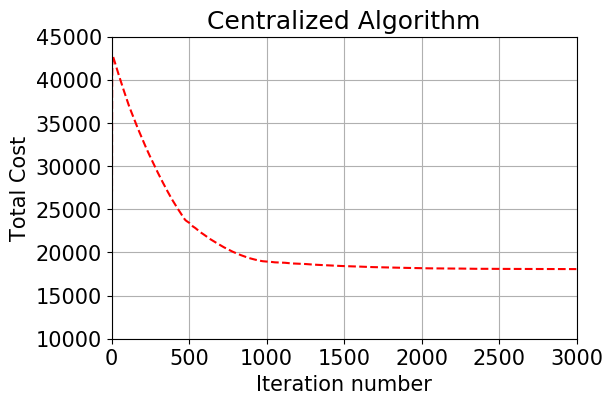}
	\includegraphics[scale=0.375]{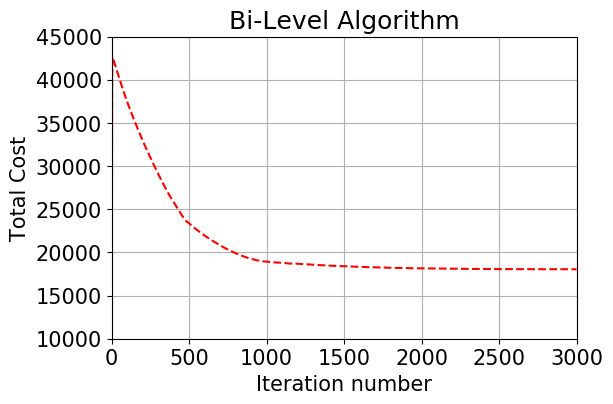}
	\includegraphics[scale=0.375]{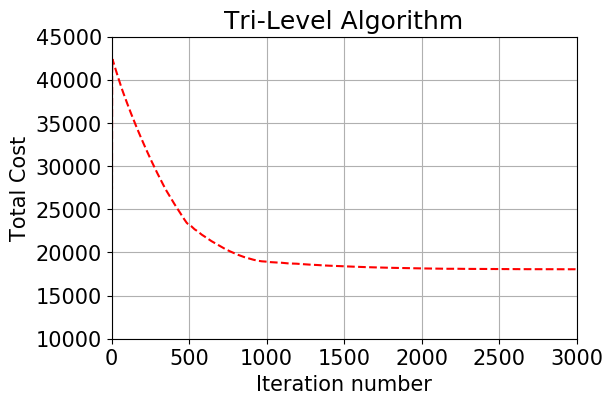}
	\caption{The centrelized algorithm (one-level algorithm), the bi-level algorithm, and the tri-level algorithms, when set to the same initialization and stepsize, exhibit identical convergence dynamics over 3,000 iterations to the same minimum cost. This is expected because the three implementations are mathematically identical.}\label{fig:converge}
\end{figure*}

\begin{table*}
\small
	\centering
	\begin{tabular}{|c|c|c|c|c|c|c|c|c|c|c|c|c|c|c|c|}
		\hline
		\multirow{2}{*}{\shortstack{One-Level\\ Algorithm }} 
		&Device \# & \multicolumn{14}{|c|}{1043}\\
		\cline{2-16}
		&Time (s) & \multicolumn{14}{|c|}{\textbf{11,588.72}} \\
		\hline
		\multirow{3}{*}{\shortstack{Bi-Level\\ Algorithm}} &Areas & \multicolumn{4}{|c|}{Area 1} & \multicolumn{4}{|c|}{Area 2} & \multicolumn{4}{|c|}{Area 3} & Area 4& Total\\
		\cline{2-16}
		&Device \# & \multicolumn{4}{|c|}{357} & \multicolumn{4}{|c|}{222} & \multicolumn{4}{|c|}{310} & 154 & 1043\\
		\cline{2-16}
		&Time (s) & \multicolumn{4}{|c|}{\textbf{1,050.16}} & \multicolumn{4}{|c|}{\textbf{501.21}} & \multicolumn{4}{|c|}{\textbf{834.90}} & \textbf{256.81} & \textbf{2,717.25}\\
		\hline
		\multirow{4}{*}{\shortstack{Tri-Level\\ Algorithm}} & Subareas & $1_1$ & $1_2$ & $1_3$ & $1_0$ & $2_1$ & $2_2$ & $2_3$ & $2_0$ &  $3_1$ & $3_2$ & $3_3$ & $3_0$ & $4_0$ & Total \\
		\cline{2-16}
		& Device \# & 49 & 74 & 23 & 211 & 70 & 39 & 17 & 96 &  68 & 66 & 68 & 108 & 154 & 1043 \\
		\cline{2-16}
		& Time (s) & \multicolumn{4}{|c|}{\textbf{713.08}} & \multicolumn{4}{|c|}{\textbf{337.04}} &\multicolumn{4}{|c|}{\textbf{490.51}} &  \textbf{261.07} & \textbf{1,890.40} \\
		\hline
	\end{tabular}
	\caption{Detailed clustering information and computational time of all areas and subareas over 3,000 iterations for the centralized algorithm, the bi-level algorithm, and the tri-level algorithm.} 
	\label{table:result}
\end{table*}

\section{Numerical Results}\label{sec:numerical}

\begin{figure}
	\centering
	\includegraphics[trim={0 0 0 0.75cm},clip,scale=0.54]{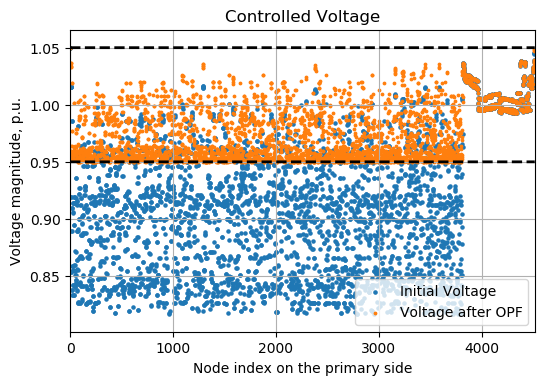}
	\caption{The largely undervoltage system has been improved after OPF.}\label{fig:voltage}
\end{figure}

\subsection{Simulation Setup}
We construct a three-phase, unbalanced, 11,000-node test feeder by connecting an IEEE 8,500-node test feeder and an EPRI Ckt7 test feeder at the substation. We simplify the network by lumping the loads on the secondary side into corresponding distribution transformers, resulting in a 4,521-node network ($N=4,521$) with a total of 1,043 controllable devices. We group all the controllable nodes into four subtree areas, and we divide Areas 1--3 further into subareas featuring subsubtrees and remaining buses; see Fig.~\ref{fig:multilevel} and Table~\ref{table:result} for the details. 

For each controllable device at bus $i$ of phase $\phi$, we consider minimizing the cost of its deviation from its initial (the most preferred) level $(p_i^{\phi}(0), q_i^{\phi}(0))$, i.e., $C^{\phi}_i(p^{\phi}_i,q^{\phi}_i)=(p^{\phi}_i-p_i^{\phi}(0))^2+(q^{\phi}_i-q_i^{\phi}(0))^2$. We set voltage bounds uniformly to $0.95$~p.u. and $1.05$~p.u. as the hard constraints. We set the same stepsizes for all three implementations as $3.5\times 10^{-4}$ for the primal updates and $3.5\times 10^{-3}$ for the dual updates. We disabled all voltage regulators so that they do not assist OPF with voltage regulation. As a result, a largely undervoltage scenario is observed, as plotted with blue dots in Fig.~\ref{fig:voltage}. We implement the primal-dual gradient algorithm to solve for the unique saddle point of the regularized Lagrangian~\eqref{eq:langr2} based on the centralized implementation, the bi-level implementation, and the tri-level implementation. OpenDSS is used to generate the nonlinear power flow at each iteration to replace the linear power flow in Eqs.~\eqref{eq:primaldual3}. 

The simulation runs with Python 3.6 on Windows 10 Enterprise Version on a laptop with Intel Core i7-7600U CPU @ 2.80GHz 2.90GHz, 8.00GB RAM. 

\subsection{Simulation Results}
\subsubsection{Dynamic Performance}
As shown in Fig.~\ref{fig:converge}, when set with the identical initial conditions and stepsizes, the centralized implementation, i.e., the straightforward implementation of Eqs.~\eqref{eq:primaldual3}, the bi-level and the tri-level implementations of Eqs.~\eqref{eq:primaldual3} share exactly the same convergence dynamics and the minimum cost value at convergence, illustrating the mathematical equivalence among the three implementation. 

\subsubsection{Speed Improvement}
Because of the computational complexity reduction achieved by our design, however, each iteration costs less time in the bi-level or the tri-level implementation than the centralized one. We record the time used to run 3,000 iterations for the three implementations in TABLE~\ref{table:result}. As shown, the time used for the centralized algorithm is 4.26 times of the time for the bi-level algorithm, which is 44\% more than the time used for the tri-lvel algorithm. 
Note that Area~4 serving as a control group has a similar computational time for bi-level and tri-level algorithms  because we do not divide it into subareas. 

Such computational complexity reduction without compromising optimality gives us ``free" speed improvement. Moreover, the multi-level structure also allows for more parallel and autonomous implementation among areas and subareas and potentially further speeds up the convergence.

\subsubsection{Voltage Regulation}
As a result, we see that the voltage profiles are regulated to within the prescribed upper and lower bounds after the OPF is solved. Note that quite a few nodes have their voltage on the lower bound 0.95~p.u. This is because intuitively the most optimal solutions should satisfy the hard voltage constraints by incurring the minimal power setpoint deviation from the initial values for devices.

\section{Conclusion}\label{sec:conclusion}
By exploring the fractal properties of power distribution networks of tree topology, as well as the structure of the impedance matrix, this paper extends our previous bi-level distributed OPF solver based on the primal-dual gradient algorithm to multi-level implementation. Privacy, including nodal voltages and network topology, has been preserved within areas and subareas. Numerical results on a 4,521-node large multi-phase distribution network show significant computational speed improvement of such an extension without compromising optimality. 


\bibliographystyle{IEEEtran}
\bibliography{biblio.bib}

\begin{thebibliography}{10}
\providecommand{\url}[1]{#1}
\csname url@samestyle\endcsname
\providecommand{\newblock}{\relax}
\providecommand{\bibinfo}[2]{#2}
\providecommand{\BIBentrySTDinterwordspacing}{\spaceskip=0pt\relax}
\providecommand{\BIBentryALTinterwordstretchfactor}{4}
\providecommand{\BIBentryALTinterwordspacing}{\spaceskip=\fontdimen2\font plus
\BIBentryALTinterwordstretchfactor\fontdimen3\font minus
  \fontdimen4\font\relax}
\providecommand{\BIBforeignlanguage}[2]{{%
\expandafter\ifx\csname l@#1\endcsname\relax
\typeout{** WARNING: IEEEtran.bst: No hyphenation pattern has been}%
\typeout{** loaded for the language `#1'. Using the pattern for}%
\typeout{** the default language instead.}%
\else
\language=\csname l@#1\endcsname
\fi
#2}}
\providecommand{\BIBdecl}{\relax}
\BIBdecl

\bibitem{zhou2019accelerated}
X.~Zhou, Z.~Liu, C.~Zhao, and L.~Chen, ``Accelerated voltage regulation in
  multi-phase distribution networks based on hierarchical distributed
  algorithm,'' \emph{IEEE Trans. Power Systems}, 2019.

\bibitem{dall2018optimala}
E.~Dall’Anese and A.~Simonetto, ``Optimal power flow pursuit,'' \emph{IEEE
  Trans. on Smart Grid}, vol.~9, no.~2, pp. 942--952, 2018.

\bibitem{tang2017real}
Y.~Tang, K.~Dvijotham, and S.~Low, ``Real-time optimal power flow,'' \emph{IEEE
  Trans. on Smart Grid}, vol.~8, no.~6, pp. 2963--2973, 2017.

\bibitem{hauswirth2017online}
A.~Hauswirth, A.~Zanardi, S.~Bolognani, F.~D{\"o}rfler, and G.~Hug, ``Online
  optimization in closed loop on the power flow manifold,'' pp. 1--6, 2017.

\bibitem{zhou2017discrete}
X.~Zhou, E.~Dall’Anese, and L.~Chen, ``Online stochastic optimization of
  networked distributed energy resources,'' \emph{IEEE Trans. on Automatic
  Control}, vol.~65, no.~6, pp. 2387--2401, 2019.

\bibitem{peng2016distributed}
Q.~Peng and S.~H. Low, ``Distributed optimal power flow algorithm for radial
  networks, i: Balanced single phase case,'' \emph{IEEE Trans. on Smart Grid},
  vol.~9, no.~1, pp. 111--121, 2016.

\bibitem{vsulc2014optimal}
P.~{\v{S}}ulc, S.~Backhaus, and M.~Chertkov, ``Optimal distributed control of
  reactive power via the alternating direction method of multipliers,''
  \emph{IEEE Transactions on Energy Conversion}, vol.~29, no.~4, pp. 968--977,
  2014.

\bibitem{magnusson2019voltage}
S.~Magn{\'u}sson, G.~Qu, C.~Fischione, and N.~Li, ``Voltage control using
  limited communication,'' \emph{IEEE Trans. on Control of Network Systems},
  vol.~6, no.~3, pp. 993--1003, 2019.

\bibitem{wu2018smart}
C.~Wu, G.~Hug, and S.~Kar, ``Smart inverter for voltage regulation: Physical
  and market implementation,'' \emph{IEEE Trans. on Power Systems}, vol.~33,
  no.~6, pp. 6181--6192, 2018.

\bibitem{kraning2014dynamic}
M.~Kraning, E.~Chu, J.~Lavaei, S.~P. Boyd \emph{et~al.}, \emph{Dynamic network
  energy management via proximal message passing}.\hskip 1em plus 0.5em minus
  0.4em\relax Now Publishers, 2014.

\bibitem{molzahn2017survey}
D.~K. Molzahn, F.~D{\"o}rfler, H.~Sandberg, S.~H. Low, S.~Chakrabarti,
  R.~Baldick, and J.~Lavaei, ``A survey of distributed optimization and control
  algorithms for electric power systems,'' \emph{IEEE Trans. on Smart Grid},
  vol.~8, no.~6, pp. 2941--2962, 2017.

\bibitem{kroposki2018autonomous}
B.~Kroposki, E.~Dall'Anese, A.~Bernstein, Y.~Zhang, and B.-M. Hodge,
  ``Autonomous energy grids,'' \emph{Proc. of Hawaii International Conference
  on System Sciences}, 2018.

\bibitem{conejo1998multi}
A.~J. Conejo and J.~A. Aguado, ``Multi-area coordinated decentralized dc
  optimal power flow,'' \emph{IEEE Trans. on Power Systems}, vol.~13, no.~4,
  pp. 1272--1278, 1998.

\bibitem{nogales2003decomposition}
F.~J. Nogales, F.~J. Prieto, and A.~J. Conejo, ``A decomposition methodology
  applied to the multi-area optimal power flow problem,'' \emph{Annals of
  operations research}, vol. 120, no. 1-4, pp. 99--116, 2003.

\bibitem{lai2014decentralized}
X.~Lai, L.~Xie, Q.~Xia, H.~Zhong, and C.~Kang, ``Decentralized multi-area
  economic dispatch via dynamic multiplier-based lagrangian relaxation,''
  \emph{IEEE Trans. on Power Systems}, vol.~30, no.~6, pp. 3225--3233, 2014.

\bibitem{zhang2018dynamic}
K.~Zhang, W.~Shi, H.~Zhu, E.~Dall'Anese, and T.~Basar, ``Dynamic power
  distribution system management with a locally connected communication
  network,'' \emph{IEEE Journal of Selected Topics in Signal Processing},
  vol.~12, no.~4, pp. 673--687, 2018.

\bibitem{chang2019saddle}
C.-Y. Chang, M.~Colombino, J.~Cort{\'e}s, and E.~Dall’Anese, ``Saddle-flow
  dynamics for distributed feedback-based optimization,'' \emph{IEEE Control
  Systems Letters}, 2019.

\bibitem{gan2014convex}
L.~Gan and S.~H. Low, ``Convex relaxations and linear approximation for optimal
  power flow in multiphase radial networks,'' \emph{Proc. of Power Systems
  Computation Conference (PSCC)}, pp. 1--9, 2014.

\bibitem{gan2016online}
------, ``An online gradient algorithm for optimal power flow on radial
  networks,'' \emph{IEEE Journal on Selected Areas in Communications}, vol.~34,
  no.~3, pp. 625--638, 2016.

\bibitem{koshal2011multiuser}
J.~Koshal, A.~Nedi{\'c}, and U.~V. Shanbhag, ``Multiuser optimization:
  Distributed algorithms and error analysis,'' \emph{SIAM Journal on
  Optimization}, vol.~21, no.~3, pp. 1046--1081, 2011.

\bibitem{bertsekas1989parallel}
D.~P. Bertsekas and J.~N. Tsitsiklis, \emph{Parallel and distributed
  computation: numerical methods}.\hskip 1em plus 0.5em minus 0.4em\relax
  Prentice hall Englewood Cliffs, NJ, 1989, vol.~23.

\end{thebibliography}

\end{document}